\documentclass[10pt,twoside,reqno]{amsart}
\usepackage{amssymb}
\usepackage{multirow}
\calclayout

\numberwithin{equation}{section}
\makeatletter

\renewcommand{\@secnumfont}{\bfseries}

\renewcommand{\section}{\@startsection{section}{1}%
  {0mm}{.7\linespacing\@plus\linespacing}{.5\linespacing}
  {\normalfont\bfseries\centering}}

\newcommand{\bibsection}{\@startsection{section}{1}%
  {0mm}{.7\linespacing\@plus\linespacing}{.5\linespacing}
  {\normalfont\scshape\centering}}

\renewcommand{\@biblabel}[1]{#1.}

\newtheorem{thm}{\bf Theorem}[section]

\newtheorem{cor}[thm]{\bf Corollary}

\newtheorem{prop}{\bf Proposition
}

\begin{document}

\vspace{1.3cm}

\title[A note on type 2 Changhee and Daehee polynomials] {A note on type 2 Changhee and Daehee polynomials}

\author{Dae San Kim}
\address{Department of Mathematics, Sogang University, Seoul 121-742, Republic of Korea}
\email{dskim@sogang.ac.kr}

\author{Taekyun Kim}
\address{Department of Mathematics, Kwangwoon University, Seoul 139-701, Republic of Korea}
\email{tkkim@kw.ac.kr}

\subjclass[2010]{11B83; 11S80}
\keywords{type 2 Euler polynomials, type 2 Bernoulli polynomials, type 2 Changhee polynomials, type 2 Daehee polynomials}
\begin{abstract}
In recent years, many authors have studied Changhee and Daehee polynomials in connection with many special numbers and polynomials.
In this paper, we investigate type 2 Changhee and Daehee numbers and polynomials and give some identities for these numbers and polynomials in relation to type 2 Euler and Bernoulli numbers and polynomials. In addition, we express the central factorial numbers of the second kind in terms of type 2  Bernoulli, type 2 Changhee and type 2 Daehee numbers of negative integral orders.
\end{abstract}

\maketitle

\section{Introduction}

Let $p$ be a fixed prime number with $p \equiv 1(\textnormal{mod}~2)$. Throughout this paper, $\mathbb{Z}_p$, $\mathbb{Q}_p$ and $\mathbb{C}_p$ will denote the ring of $p$-adic integers, the field of $p$-adic rational numbers and the completion of the algebraic closure of $\mathbb{Q}_p$, respectively. As is known, the Stirling numbers of the second kind are defined by
\begin{equation} \label{01}
\frac{1}{k!} (e^t -1)^k = \sum_{n=k}^\infty S_2 (n,k) \frac{t^n}{n!},~(k \geq 0),\,\,\,(\textnormal{see}\,\,[1,2,4,15]).
\end{equation}

For $n\geq 0$, the central factorial is defined as (\textnormal{see}\,\,[14])
\begin{equation} \label{02}
x^{[0]} = 1, x^{[n]} = x(x+\frac{n}{2}-1) (x+\frac{n}{2}-2) \cdots (x-\frac{n}{2}+1),\,\, (n \geq 1).
\end{equation}

The central factorial numbers of the second kind are given by
\begin{equation}  \label{03}
x^n = \sum_{k=0}^n T(n,k) x^{[k]}, ~ (n \geq 0),\,\,
(\textnormal{see}\,\,[1,2,15,17]).
\end{equation}

Let $f(x)$ be a continuous function on $\mathbb{Z}_p$. Then the fermionic $p$-adic integral on $\mathbb{Z}_p$ is defined by Kim as
\begin{equation} \label{04}
\int_{\mathbb{Z}_p} f(x)d\mu_{-1}(x)= \lim_{N\rightarrow \infty} \sum_{x=0}^{p^N-1}f(x)(-1)^x,\,\,\,
(\textnormal{see}\,\,[3,5-8]). \\
\end{equation}

Thus, by \eqref{04}, we get
\begin{equation} \label{05}
\int_{\mathbb{Z}_p} f(x+1)d\mu_{-1}(x) + \int_{\mathbb{Z}_p} f(x)d\mu_{-1}(x) = 2 f(0),\,\,\,(\textnormal{see}\,\,[5,6]). \\
\end{equation}

For $t \in \mathbb{C}_p$ with $|t|_p < p^{-\frac{1}{p-1}}$, the generating  function of Changhee polynomials can be expressed in terms of the following fermionic integral on $\mathbb{Z}_p$:
\begin{equation} \label{06}
\int_{\mathbb{Z}_p} (1+t)^{x+y} du_{-1}(y)= \frac{2}{2+t} (1+t)^x = \sum_{n=0}^\infty Ch_n(x) \frac{t^n}{n!}.
\end{equation}

In particular, for $x=0$, $Ch_n = Ch_n(0)$ are called the Changhee numbers(\textnormal{see}\,\,[7,11-14]).

From \eqref{06}, we note that
\begin{equation} \label{07}
E_n(x) = \sum_{l=0}^n S_2(n,l) Ch_l(x),~ (n \geq 0),\,\,\,(\textnormal{see}\,\,[9,12,13,16]),
\end{equation}

\noindent where $E_n(x)$ are the ordinary Euler polynomials given by
\begin{equation} \label{08}
\frac{2}{e^t+1} e^{xt} = \sum_{n=0}^\infty E_n(x) \frac{t^n}{n!},\,\,\,(\textnormal{see}\,\,[4-10]).
\end{equation}

From \eqref{08}, we note that
\begin{equation}
Ch_n(x) = \sum_{l=0}^n E_l(x) S_1(n,l),\,\,\,(\textnormal{see}\,\,[9,12,16]),
\end{equation} \label{09}
where $S_1(n,l)$ are the Stirling numbers of the first kind.

In this paper, we introduce type 2 Changhee numbers and polynomials which can be expressed in terms of  fermionic $p$-adic integrals on $\mathbb{Z}_p$.
We derive some identities for these numbers and the polynomials in relation to type 2 Euler and Bernoulli numbers and polynomials. In addition, we express the central factorial numbers of the second kind in terms of type 2  Bernoulli, type 2 Changhee and type 2 Daehee numbers of negative integral orders.

\vspace{0.1in}

\section{Type 2 Changhee and Daehee numbers and polynomials}
\vspace{0.1in}

In this section, we assume that $t \in \mathbb{C}_p$ with $|t|_p < p^{-\frac{1}{p-1}}$.
As is well known, the type 2 Euler polynomials are defined by the generating function
\begin{equation}\label{10}
\int_{\mathbb{Z}_p} e^{(2y+1+x)t} d\mu_{-1}(y)= \frac{2}{e^t + e^{-t}} e^{xt} = \sum_{n=0}^\infty E_n^{*}(x) \frac{t^n}{n!}.
\end{equation}

Indeed, we note that $E_n^{*}(x) = 2^nE_n(\frac{x+1}{2}), ~(n \geq 0),\,\,\,(\textnormal{see}\,\,[4-10])$.
When $x=0$, $E_n^{*} = E_n(0)$ will be called the type 2 Euler numbers in this paper.

Motivated by \eqref{06} and \eqref{10}, we define the type 2 Changhee polynomials as
\begin{equation} \label{11}
\int_{\mathbb{Z}_p} (1+t)^{2y+1+x} du_{-1}(y)= \frac{2}{(1+t)+(1+t)^{-1}} (1+t)^x = \sum_{n=0}^\infty c_n(x) \frac{t^n}{n!}.
\end{equation}
When $x=0$, $c_n = c_n(0)$ are called the type 2 Changhee numbers.

By replacing $t$ by $\log(1+t)$ in \eqref{10}, we get
\begin{equation}\label{12}
\begin{split}
\frac{2}{(1+t)+(1+t)^{-1}} (1+t)^x  & = \sum_{l=0}^\infty E_l^{*}(x) \frac{1}{l!} (\log(1+t))^l \\
                                    & = \sum_{n=0}^\infty \Big( \sum_{l=0}^n E_l^{*}(x) S_1(n,l) \Big) \frac{t^n}{n!},
\end{split}
\end{equation}
where $S_1(n,l)$ are the Stirling numbers of the first kind.

From \eqref{11} and \eqref{12}, we have the following theorem.
\begin{thm}
For $n \geq 0$, we have
\begin{equation*}\label{Thm1}
c_n(x) =  \sum_{l=0}^n E_l^{*}(x) S_1(n,l).
\end{equation*}
\end{thm}

\vspace{0.1in}

Replacing $t$ by $e^t -1$ in \eqref{11}, we have
\begin{equation} \label{13}
\frac{2}{e^t + e^{-t}} e^{xt} = \sum_{l=0}^\infty c_l(x) \frac{1}{l!} (e^t -1)^l =  \sum_{n=0}^\infty \Big( \sum_{l=0}^n S_2(n,l) c_l(x) \Big) \frac{t^n}{n!}.
\end{equation}

Therefore, by \eqref{10} and \eqref{13}, we get the following theorem.
\begin{thm}
For $n \geq 0$, we have
\begin{equation*}\label{Thm2}
E_n^{*}(x) = \sum_{l=0}^n S_2(n,l) c_l(x).
\end{equation*}
\end{thm}

\vspace{0.1in}

For $\alpha \in \mathbb{R}$, let us define the type 2 Changhee polynomials of order $\alpha$ by
\begin{equation} \label{14}
\Big(\frac{2}{(1+t)+(1+t)^{-1}}\Big)^{\alpha} (1+t)^x = \sum_{n=0}^\infty c_n^{(\alpha)}(x) \frac{t^n}{n!}.
\end{equation}
When $x=0$, $C_n^{(\alpha)} = C_n^{(\alpha)}(0)$ are called the type 2 Changhee numbers of order $\alpha$.

For $k \in \mathbb{N} \cup \{0\}$, let us take $\alpha = -k$. Then, by \eqref{14} with $x=0$ and $t$ replaced by $e^{\frac{t}{2}}-1$, we get
\begin{equation}\label{15}
\begin{split}
\frac{1}{2^k} (e^{\frac{t}{2}} - e^{-\frac{t}{2}})^k & = \sum_{l=0}^\infty c_l^{(-k)} \frac{1}{l!} (e^{\frac{t}{2}} -1)^l \\
& = \sum_{n=0}^\infty \Big( \frac{1}{2^n} \sum_{l=0}^n S_2(n,l) c_l^{(-k)} \Big) \frac{t^n}{n!}.
\end{split}
\end{equation}

On the other hand, by Proposition 1 to be shown below, we get
\begin{equation}\label{16}
\frac{1}{2^k}(e^{\frac{t}{2}}-e^{-\frac{t}{2}})^k=\frac{k!}{2^k} \frac{1}{k!} (e^{\frac{t}{2}} - e^{-\frac{t}{2}})^k  = \frac{k!}{2^k} \sum_{n=k}^\infty T(n,k) \frac{t^n}{n!}.
\end{equation}

Therefore, by \eqref{15} and \eqref{16}, we get the following theorem.
\begin{thm}
For $n \geq k$, we have
\begin{equation*}\label{Thm3}
T(n,k)  = \frac{2^{k-n}}{k!} \sum_{l=0}^n S_2(n,l) c_l^{(-k)},
\end{equation*}
and, for $0 \leq n \leq k-1$,
\begin{equation*}
\sum_{l=0}^n S_2(n,l) c_l^{(-k)}=0.
\end{equation*}
\end{thm}

\vspace{0.1in}

The central difference operator is given by (\textnormal{see}\,\,[14,15])
\begin{equation}\label{17}
\begin{split}
& \delta f(x) = f(x+\frac{1}{2}) - f(x-\frac{1}{2}).
\end{split}\end{equation}

\noindent After appying the operator $\delta$ once more, we obtain
\begin{equation}\begin{split}
& \delta^2 f(x) = f(x+1) - 2 f(x) + f(x-1) = \sum_{l=0}^2 \binom{2}{l} (-1)^{2-l} f(x+l- \frac{2}{2}).
\end{split}
\end{equation}

\noindent Continuing this process, we have
\begin{equation}\label{18}
\delta^k f(x) = \sum_{l=0}^k \binom{k}{l} (-1)^{k-l} f(x+l- \frac{k}{2}), ~(k \in \mathbb{N} \cup \{0\}).
\end{equation}

Let $f(x)$ be analytic at $x=b, ~( b \in \mathbb{R})$.
Then $f(x)$ can be rewritten as
\begin{equation}\label{19}
f(x) = \sum_{n=0}^\infty  A_n (x-b)^{[n]}.
\end{equation}

Now, we observe that
\begin{equation}\label{20}
\begin{split}
\delta x^{[n]} & = (x+\frac{1}{2})^{[n]} - (x-\frac{1}{2})^{[n]} \\
               & = (x+\frac{1}{2}) (x+\frac{n-1}{2})(x+\frac{n-3}{2}) \cdots (x-\frac{n-3}{2}) - (x-\frac{1}{2}) (x+\frac{n-3}{2}) \cdots (x-\frac{n-1}{2}) \\
               & = (x+\frac{n-3}{2}) (x+\frac{n-4}{2}) \cdots (x-\frac{n-3}{2}) \Big\{(x+\frac{1}{2}) (x+\frac{n-1}{2}) - (x-\frac{1}{2}) (x-\frac{n-1}{2}) \Big\} \\
               & = n x (x+\frac{n-3}{2}) (x+\frac{n-4}{2}) \cdots (x-\frac{n-3}{2}) \\
               & = n x (x+\frac{n-1}{2}-1) (x+\frac{n-1}{2}-2) \cdots (x-\frac{n-1}{2}+1) \\
               & = n x^{[n-1]}.
\end{split}
\end{equation}

From \eqref{19} and \eqref{20}, we have
\begin{equation}\label{21}
\delta^k f(b) = \delta^k f(x)|_{x=b} =  A_k k!.
\end{equation}

By \eqref{19} and \eqref{20}, we get
\begin{equation}\label{22}
f(x) = \sum_{k=0}^\infty  \frac{\delta^k f(b)}{k!} (x-b)^{[k]}.
\end{equation}

In particular, if we take $b=0$, then we have
\begin{equation}\label{23}
f(x) = \sum_{k=0}^\infty  \frac{\delta^kf(0)}{k!} x^{[k]}.
\end{equation}

Let us take $f(x) = e^{xt}$. Then, by \eqref{18}, we get
\begin{equation}\label{24}
\begin{split}
\delta^k f(0) & = \sum_{l=0}^k \binom{k}{l}(-1)^{k-l} e^{(l-\frac{k}{2}) t} \\
              & = e^{-\frac{k}{2} t} (e^t -1)^k = ( e^{\frac{t}{2}} - e^{-\frac{t}{2}})^k.
\end{split}
\end{equation}

From \eqref{23} and \eqref{24}, we have
\begin{equation}\label{25}
e^{xt} = \sum_{k=0}^\infty  \frac{1}{k!} ( e^{\frac{t}{2}} - e^{-\frac{t}{2}})^k x^{[k]}.
\end{equation}

On the other hand, by \eqref{03}, we get
\begin{equation}\label{26}
\begin{split}
e^{xt} & = \sum_{n=0}^\infty x^n \frac{t^n}{n!} = \sum_{n=0}^\infty \Big( \sum_{k=0}^n T(n,k) x^{[k]} \Big)  \frac{t^n}{n!}  \\
       & =  \sum_{k=0}^\infty \Big( \sum_{n=k}^\infty T(n,k) \frac{t^n}{n!}  \Big) x^{[k]}.
\end{split}
\end{equation}

Therefore, by \eqref{25} and \eqref{26}, we obtain the following proposition.
\begin{prop}
For $k \geq 0$, we have
\begin{equation*}\label{prop1}
\frac{1}{k!} ( e^{\frac{t}{2}} - e^{-\frac{t}{2}})^k  = \sum_{n=k}^\infty T(n,k) \frac{t^n}{n!}.
\end{equation*}
\end{prop}

It is well known that Daehee polynomials are defined by
\begin{equation} \label{27}
\frac{\log(1+t)}{t} (1+t)^x = \sum_{n=0}^\infty D_n(x) \frac{t^n}{n!},\,\,\,(\textnormal{see}\,\,[4,10]).
\end{equation}
When $x=0$, $D_n = D_n(0)$ are called the Daehee numbers.

Now, we consider the type 2 Bernoulli polynomials given by
\begin{equation}\label{28}
\frac{t}{e^t - e^{-t}} e^{xt} = \sum_{n=0}^\infty b_n (x) \frac{t^n}{n!}.
\end{equation}
When $x=0$, $b_n = b_n(0)$ are called the type 2 Bernoulli numbers.

Indeed, we note that
\begin{equation}\label{29}
b_n(x) = 2^{n-1} B_n(\frac{x+1}{2}), ~ (n \geq 0),
\end{equation}
where $B_n(x)$ are the ordinary Bernoulli polynomials defined by
\begin{equation} \label{30}
\frac{t}{e^t -1} e^{xt} = \sum_{n=0}^\infty B_n(x) \frac{t^n}{n!},\,\,\,(\textnormal{see}\,\,[5,15]).
\end{equation}

\vspace{0.1in}

In the view of \eqref{27} and \eqref{28}, we define the type 2 Daehee polynomials by
\begin{equation} \label{31}
\frac{\log(1+t)}{(1+t)-(1+t)^{-1}} (1+t)^x = \sum_{n=0}^\infty d_n(x) \frac{t^n}{n!}.
\end{equation}
When $x=0$, $d_n = d_n(0)$ are called the type 2 Daehee numbers.

Replacing $t$ by $e^t -1$ in \eqref{31}, we get
\begin{equation}\label{32}
\begin{split}
\frac{t}{e^t - e^{-t}} e^{xt} & = \sum_{l=0}^\infty d_l(x) \frac{1}{l!} (e^t -1)^l \\
& = \sum_{n=0}^\infty \Big( \sum_{l=0}^n S_2(n,l) d_l(x) \Big) \frac{t^n}{n!},
\end{split}
\end{equation}
where $S_2(n,l)$ are the Stirling numbers of the second kind.

From \eqref{28} and \eqref{32}, we have the following theorem.
\begin{thm}
For $n \geq 0$, we have
\begin{equation*}\label{Thm10}
 b_n(x) =  \sum_{l=0}^n S_2(n,l) d_l(x).
\end{equation*}
In particular,
\begin{equation*}
 b_n =  \sum_{l=0}^n S_2(n,l) d_l.
\end{equation*}
\end{thm}

\vspace{0.1in}

From \eqref{28}, we can derive the following equation.
\begin{equation}\label{33}
\begin{split}
\frac{\log(1+t)}{(1+t)-(1+t)^{-1}} (1+t)^x & = \sum_{l=0}^\infty b_l(x) \frac{1}{l!} (\log(1+t))^l \\
                                    & = \sum_{n=0}^\infty \Big(\sum_{l=0}^n S_1(n,l) b_l(x) \Big) \frac{t^n}{n!},
\end{split}
\end{equation}

From \eqref{31} and \eqref{33}, we obtain the following theorem.
\begin{thm}
For $n \geq 0$, we have
\begin{equation*}\label{Thm11}
 d_n(x) =  \sum_{l=0}^n S_1(n,l) b_l(x).
\end{equation*}
In particular,
\begin{equation*}
 d_n =  \sum_{l=0}^n S_1(n,l) b_l.
\end{equation*}
\end{thm}

\vspace{0.1in}

For $\alpha \in \mathbb{R}$, let us define the type 2 Daehee polynomials of order $\alpha$ by
\begin{equation} \label{34}
\Big( \frac{\log(1+t)}{(1+t)-(1+t)^{-1}}\Big)^{\alpha} (1+t)^x = \sum_{n=0}^\infty d_n^{(\alpha)}(x) \frac{t^n}{n!}.
\end{equation}
When $x=0$, $d_n^{(\alpha)} = d_n^{(\alpha)}(0)$ are called the type 2 Daehee numbers of order $\alpha$.

For $k \in \mathbb{N} \cup \{0\}$, let us take $\alpha = -k$. Then, by \eqref{34} with $x=0$ and $t$ replaced by $e^{\frac{t}{2}}-1$, we get
\begin{equation}\label{35}
\begin{split}
(\frac{2}{t})^k (e^{\frac{t}{2}} - e^{-\frac{t}{2}})^k  & = \sum_{l=0}^\infty d_l^{(-k)} \frac{1}{l!} (e^{\frac{t}{2}} -1)^l \\
                                                        & = \sum_{n=0}^\infty \Big( \frac{1}{2^n} \sum_{l=0}^n d_l^{(-k)} S_2(n,l) \Big) \frac{t^n}{n!}.
\end{split}
\end{equation}

On the other hand, by Proposition 1, we get
\begin{equation}\label{36}
\begin{split}
(\frac{2}{t})^k (e^{\frac{t}{2}} - e^{-\frac{t}{2}})^k  & = \frac{2^k k!}{t^k} \frac{1}{k!} (e^{\frac{t}{2}} - e^{-\frac{t}{2}})^k \\
                                                        & = \frac{2^k k!}{t^k} \sum_{n=k}^\infty T(n,k) \frac{t^n}{n!} \\
                                                        & = 2^k k! \sum_{n=0}^\infty T(n+k,k) \frac{t^n}{(n+k)!} \\
                                                        & = 2^k \sum_{n=0}^\infty T(n+k,k) \frac{1}{\binom{n+k}{k}} \frac{t^n}{n!} \\
\end{split}
\end{equation}

Therefore, by \eqref{35} and \eqref{36}, we get the following theorem.
\begin{thm}
For $n,k \geq 0$, we have
\begin{equation*}\label{Thm12}
2^{n+k} T(n+k,k) =  \binom{n+k}{k} \sum_{l=0}^n d_l^{(-k)} S_2(n,l).
\end{equation*}
\end{thm}

\vspace{0.1in}

Now, we define the type 2 Bernoulli polynomials of order $\alpha$ by
\begin{equation}\label{37}
\Big(\frac{t}{e^t - e^{-t}}\Big)^{\alpha} e^{xt} = \sum_{n=0}^\infty b_n^{(\alpha)} (x) \frac{t^n}{n!}.
\end{equation}

In particular, $\alpha = k \in \mathbb{N}$, we have
\begin{equation*}
t^k e^{xt} \underbrace{\text{csch} t \times \cdots \times \text{csch} t}_{k~\text{times}} = 2^k \sum_{n=0}^\infty b_n^{(\alpha)} (x) \frac{t^n}{n!}.
\end{equation*}
When $x=0$, $b_n^{(\alpha)} = b_n^{(\alpha)}(0)$ are called the type 2 Bernoulli numbers of order $\alpha$.

From by \eqref{34} and \eqref{37}, we have the following Corollary.
\begin{cor}
For $n \geq 0$, we have
\begin{equation*}\label{cor1}
b_n^{(\alpha)} (x) =  \sum_{l=0}^n S_2(n,l) d_l^{(\alpha)}(x),
\end{equation*}
and
\begin{equation*}
d_n^{(\alpha)}(x) =  \sum_{l=0}^n S_1(n,l) b_l^{(\alpha)}(x).
\end{equation*}
\end{cor}

\vspace{0.1in}

From \eqref{37}, we observe that
\begin{equation}\label{38}
\begin{split}
\frac{1}{t^k} (e^{\frac{t}{2}} - e^{-\frac{t}{2}})^k  =  \sum_{n=0}^\infty \frac {b_n^{(-k)}}{2^{n+k}} \frac{t^n}{n!}. \\
\end{split}
\end{equation}

On the other hand,
\begin{equation}\label{39}
\begin{split}
\frac{1}{t^k} (e^{\frac{t}{2}} - e^{-\frac{t}{2}})^k  & = \frac{k!}{t^k} \sum_{n=k}^\infty  T(n,k) \frac{t^n}{n!} \\
& =  \sum_{n=0}^\infty \frac{T(n+k,k)}{\binom{n+k}{k}} \frac{t^n}{n!}. \\
\end{split}
\end{equation}

Therefore, by \eqref{38} and \eqref{39}, we get the following theorem.
\begin{thm}
For $n,k \geq 0$, we have
\begin{equation*}\label{Thm13}
2^{n+k} T(n+k,k) =  \binom{n+k}{k} b_n^{(-k)}.
\end{equation*}
\end{thm}


\vspace{0.1in}

\begin{thebibliography}{0}

\bibitem{01} L. Comtet,
\textit{Nombres de Stirling g\'en\'eraux et fonctions sym\'etriques,
} C. R. Acad, Sci. Paris Ser. A., {\bf{275}} (1972), 747--750.

\bibitem{02} L. Comtet,
\textit{Advanced combinatorics: the art of finite and infinite expansions( translated from the French by J.W. Nienhuys),
} Dordrecht and Boston: Reidel, {\bf{1974}}.

\bibitem{03} D. V. Dolgy, G. -W. Jang, D. S. Kim, T. Kim,
\textit{Explicit expressions for Catalan-Daehee numbers,
} Proc. Jangjeon Math. Soc., {\bf{20}} (2017), no. 1, 1--9.

\bibitem{04} B. S. El-Desouky, A. Mustafa,
\textit{New results on higher-order Daehee and Bernoulli numbers and polynomials,
} Adv. Difference Equ., 2016, paper No. 32, 21 pp.

\bibitem{05} D. S. Kim, T. Kim,
\textit{Some $p$-adic integral on $\mathbb{Z}_p$ associated with trigonometric functions,
} Russ. J. Math. Phys., {\bf{25}} (2018), no. 3, 300--308.

\bibitem{06} D. S. Kim, T. Kim, H. I. Kwon, G. -W. Jang,
\textit{Degenerate Daehee polynomials of the second kind,
} Proc. Jangjeon Math. Soc., {\bf{21}} (2018), no. 1, 83--97.

\bibitem{07} T. Kim,
\textit{Lebesgue-Radon-Nikodym theorem with respect to fermionic $p$-adic invariant measure on $\mathbb{Z}_p$,
} Adv. Stud. Contemp. Math.(Kyungshang), ${\mathbf{28}}$ (2018), no. 1, 149--160.

\bibitem{08} T. Kim,
\textit{$q$-Volkenborn integration,
} Russ. J. Math. Phys., {\bf{9}} (2002), no. 3, 288--299.

\bibitem{09} T. Kim, D. S. Kim,
\textit{A note on nonlinear Changhee differential equations,
} Russ. J. Math. Phys., {\bf{23}} (2016), no. 1, 88--92.

\bibitem{10} W. A. Khan, K. S. Nisar, U. Duran, M. Acikgoz, S. Araci,
\textit{Multifarious implicit summation formulae of Hermite-based poly-Daehee polynomials,
} Appl. Math. Inf. Sci., ${\mathbf{12}}$ (2018), no. 2, 305--310.

\bibitem{11} H. -I. Kwon, T. Kim, J. J. Seo,
\textit{A note on degenerate Changhee numbers and polynomials,
} Proc. Jangjeon Math. Soc., {\bf{18}} (2015), no. 3, 295--305.

\bibitem{12} C. Liu, W. Wuyungaowa,
\textit{Application of probabilistic method on Daehee sequences,
} Eur. J. Pure Appl. Math., ${\mathbf{11}}$ (2018), no. 1, 69--78.

\bibitem{13} H. -K. Pak, J. Jeong, D. -J. Kang, S. -H. Rim,
\textit{Changhee-Genocchi numbers and their applications,
} Ars Combin., {\bf{136}} (2018), 153--159.

\bibitem{14} J. Riordan,
\textit{Combinatorial Identities,} New York: Wiley, {\bf{1968}}.

\bibitem{15} S. Roman,
\textit{The Umbral Calculus,} New York: Academic Press, {\bf{1984}}.

\bibitem{16} Y. Simsek,
\textit{Identities on the Changhee numbers and Apostol-type Daehee polynomials,
} Adv. Stud. Contemp. Math.(Kyungshang), ${\mathbf{27}}$ (2017), no. 2, 199--212.

\bibitem{17} W.P. Zhang,
\textit{Some identities involving the Euler and central factorial numbers,} Fibonacci Quart.
{\bf{36}} (1998), 154--157.


\end{thebibliography}
\end{document}